\documentclass[12pt]{article}
\usepackage[utf8]{inputenc}
\usepackage[T1]{fontenc}

\usepackage[margin=1in]{geometry}

\usepackage{pslatex}
\usepackage[compact]{titlesec}

\usepackage{amsmath, amssymb, amsthm}

\newtheorem{theorem}{Theorem}[section]
\newtheorem{lemma}{Lemma}
\newtheorem{case3surface}{Case}
\newtheorem{case2surface}{Case}
\newtheorem{casezero}{Case}

\begin{document}

\title{Line-Plane Incidence Bound in $\mathbb{R}^4$}
\date{September 25, 2017}
\author{Chao Cheng}
\maketitle

\begin{abstract}
We consider an incidence problem in $\mathbb{R}^4$ which asks, for a set of $L$ lines and a set of $S$ planes in general position, what the maximum number of line-plane incidences is. A line-plane incidence is defined as a point where a line and a plane intersect. We prove that, when the lines and planes are in a truly 4-dimensional configuration such that no more than $L^{\frac{1}{2}+\epsilon}$ lines are contained in any 2-dimensional surface of degree at most $D$ and no more than $S^{\frac{1}{2}+\epsilon}$ 2-planes are contained in any 3-dimensional hypersurface of degree at most $D$, and if $L^{1/2} \ll S \ll L$, then for a constant $D>1$ and an $\epsilon>0$ there exists a non-trivial upper bound for incidences between lines and planes: $L^{\frac{3}{4}+\frac{1}{2}\epsilon}S + LS^{\frac{1}{2}+\epsilon}$. We also prove several supporting lemmas. 
\end{abstract}

\thispagestyle{empty}

\newpage

\setcounter{page}{1}

\section{Introduction}
Incidence Geometry is the study of intersection patterns between simple geometric objects. A fundamental question to ask in this field is how various objects in space---points, lines, circles, planes---can intersect one another. In their celebrated work \cite{szemereditrotter}, Szemer\'edi and Trotter proved the following tight upper bound on the number of incidences between a set of points and a set of lines, where an incidence is defined as a point-line pair such that the point lies on the line.

\begin{theorem}[Szemer\'edi-Trotter]
Given $n$ points and $m$ lines in a plane, the number of incidences is $O(n^\frac{2}{3}m^\frac{2}{3}+n+m)$.
\end{theorem}

Due to its broad applications and significance in the field, the Szemer\'edi–Trotter theorem has led to extensive study of similar types of incidence questions. A number of other mathematical problems have been transformed to questions about Szemer\'edi-Trotter type incidence bounds \cite{solymositao}. For example, Elekes used the Szemer\'edi-Trotter theorem to obtain new bounds on the sum-product problem. In 2005, Bennett, Carbery, and Tao discovered a connection between bounds on the number of incidences between points and lines in $\mathbb{R}^3$ and multilinear Kakeya estimates. Recently, Guth and Katz also used bounds on the number of incidences between points and lines in $\mathbb{R}^3$ in their solution to the Erd\H{o}s distinct distances problem. Further applications of Szemer\'edi-Trotter type incidence bounds in mathematics and theoretical computer science, including the construction of extractors used in cryptography and data structures, are discussed in \cite{dvir}. 

Although the Szemer\'edi-Trotter theorem has inspired other related bounds, including bounds for point-line incidences in higher dimensions and for point-plane incidences, an incidence bound between lines and 2-dimensional planes has not yet been determined, which is the main objective of this project. We define an \textit{incidence} between lines and planes to be a point where a line and a plane intersect. We then consider the fundamental question of obtaining an upper bound for the incidences between lines and planes: in $\mathbb{R}^3$, intersections between lines and planes in general position will have a trivial upper bound, as each line must intersect each plane; in higher dimensions, however, a line and a plane in general position do not necessarily intersect. In this paper, we prove a non-trivial upper bound for line-plane incidences in $\mathbb{R}^4$. 

In considering incidence problems in higher dimensions, it is necessary to require that the objects cannot be projected to lower dimensions. For points and lines in $\mathbb{R}^3$, if the points and lines are not restricted to a truly 3-dimensional configuration, then it is possible to project them onto a generic 2-dimensional plane while preserving all incidence relations, giving a bound no different from the Szemer\'edi-Trotter theorem for point-line incidences in two dimensions.  In his book \cite{guthbook}, Guth was able to directly address this issue by adding conditions to ensure that 3-dimensional situations were prevented from being reduced to a 2-dimensional case, using restrictions such as "there can be no more than $B$ lines contained in any plane" for a constant $B$. Non-degeneracy conditions similar to Guth's are often used to ensure that configurations of objects are truly 3-dimensional or truly $n$-dimensional. For example, as a part of their solution to the Erd\H{o}s distinct distances problem \cite{guthkatz}, Guth and Katz included a non-degeneracy restriction: "let $\mathcal{L}$ be a set of $L$ lines in $\mathbb{R}^3$ with at most $B$ lines in any plane." Similarly, in their extension of the Szemer\'edi-Trotter theorem to points and lines in $\mathbb{R}^4$ \cite{sharirsolomon}, Sharir and Solomon included a condition that not too many lines are contained in any single plane. In this paper, we also utilize non-degeneracy restrictions to ensure 4-dimensional configurations of lines and planes in order to determine a meaningful bound; namely, we restrict the number of lines any 2-surface can contain and the number of 2-planes any 3-hypersurface can contain.

\section{Methods and Techniques}
\subsection{Polynomial Partitioning}
In their solution to the Erd\H{o}s distinct distances problem \cite{guthkatz}, Guth and Katz developed a new method of solving incidence geometry problems, known as polynomial partitioning.	

\begin{theorem}[Theorem 4.1 of \cite{guthkatz}]
If $\mathcal{A}$ is a set of $A$ points in $\mathbb{R}^n$, and $J\geq1$ is an integer, then there is a non-zero polynomial surface $Z(P)$ of degree $D \leq 2^{(J/n)}$ with the following property. The complement $\mathbb{R}^n\setminus Z(P)$ is the union of $2^J$ open cells $O_i$, and each cell contains $\leq 2^{-J}A$ points of $\mathcal{A}$.
\end{theorem}

Their technique uses a "divide and conquer" approach to counting incidences by using a polynomial $P$ of degree $D$ to partition space into a number of cells. The zero set of the polynomial, $Z(P)$, effectively divides $\mathbb{R}^n$ into at most $D^n$ cells. A key benefit of this technique is that in each cell the number of points is small, and the degree of the polynomial $P$ is not very large. In \cite{kms}, the authors used polynomial partitioning to give new, simpler proofs of several classical theorems. In addition to points, Guth later generalized the polynomial partitioning method further to other objects such as lines and 2-planes.

\subsection{B\'ezout's Theorem} 
B\'ezout's Theorem is extremely relevant to the proof. The theorem was originally stated in 1687 by Isaac Newton and said that the number of intersection points of two algebraic curves is given by the product of their degrees. B\'ezout's Theorem has since been generalized to higher dimensions for higher-dimensional objects, such as 2-planes. In \cite{guthdistinctdistances}, the $\mathbb{R}^2$ version and $\mathbb{R}^3$ version of the B\'ezout's Theorem were proven and re-stated as follows.

\begin{theorem}
If $Q_1$, $Q_2$ are non-zero polynomials in $\mathbb{R}^2$ with no common factor, then Z($Q_1$) $\cap$ Z($Q_2$) $\subset \mathbb{R}^2$ contains at most (Deg $Q_1$)(Deg $Q_2$) points.
\end{theorem}

\begin{theorem}
If $Q_1$, $Q_2$ are non-zero polynomials in $\mathbb{R}^3$ with no common factor, then Z($Q_1$) $\cap$ Z($Q_2$) $\subset \mathbb{R}^3$ contains at most (Deg $Q_1$)(Deg $Q_2$) lines.
\end{theorem}

Our proof also requires a $\mathbb{R}^4$ version of B\'ezout's Theorem, which can be derived from a generalized form of B\'ezout's Theorem for higher dimensional space in \cite{sharirsolomon}. We can restate a $\mathbb{R}^4$ version as follows.
 
\begin{theorem} (B\'ezout's Theorem for $\mathbb{R}^4$)
If $Q_1$, $Q_2$ are non-zero 3-dimensional polynomials in $\mathbb{R}^4$ with no common factor, then Z($Q_1$) $\cap$ Z($Q_2$) $\subset \mathbb{R}^4$ contains at most (Deg $Q_1$)(Deg $Q_2$) 2-planes.
\end{theorem}

\subsection{Zarankiewicz Problem}
A bipartite graph consists of two disjoint sets of vertices $\mathcal{U}$ and $\mathcal{V}$ and a set of edges that connect vertices in $\mathcal{U}$ to vertices in $\mathcal{V}$, where no two edges can connect the same pair of vertices. A complete bipartite graph, in which every vertex in $\mathcal{U}$ is connected to every vertex in $\mathcal{V}$, is denoted by $K_{s,t}$, where $s$ is the number of vertices in $\mathcal{U}$ and $t$ is the number of vertices in $\mathcal{V}$. The Zarankiewicz function $z(m, n; s, t)$ denotes the maximum possible number of edges for a bipartite graph with two sets of vertices $\mathcal{U}$ and $\mathcal{V}$ for which $|\mathcal{U}|=m$ and $|\mathcal{V}|=n$, that does not contain any complete subgraph of the form $K_{s,t}$. The Kov\'ari-S\'os-Tur\'an Theorem provides an upper bound for the Zarankiewicz function \cite{kst}.

\begin{theorem}[Kov\'ari-S\'os-Tur\'an]
$z(m, n; s,t) < (s-1)^\frac{1}{t}(n-t+1)m^{1-\frac{1}{t}} + (t-1)m$
\end{theorem}

This result is useful in situations where a bipartite graph can model the intersections between two sets of objects.

\section{Main Theorem}
We proceed to the main theorem of the paper. From here on, we refer to 2-dimensional planes as "2-planes" and 3-dimensional hypersurfaces as "3-surfaces".
\begin{theorem}[Main Theorem]
\label{main}
Given $L$ lines and $S$ 2-planes in $\mathbb{R}^4$ such that $L^{1/2} \ll S \ll L$, with at most $L^{\frac{1}{2}+\epsilon}$ lines in any irreducible algebraic 2-surface of degree at most $D$ and at most $S^{\frac{1}{2}+\epsilon}$ 2-planes in any irreducible algebraic 3-surface degree at most $D$, then for a constant $D>1$ and an $\epsilon>0$, the intersections between lines and 2-planes is bounded by $O(L^{\frac{3}{4}+\frac{1}{2}\epsilon}S + LS^{\frac{1}{2}+\epsilon})$.
\end{theorem}

We will use two techniques to prove our main theorem: polynomial partitioning and induction. In \cite{guthdistinctdistances}, Guth introduced the framework of using a slightly stronger bound to make the inductive argument close in a similar situation, where he proved a bound for $r$-rich points in $\mathbb{R}^3$. The same technique can used in our case. Our slightly stronger bound says that given $L$ lines and $S$ 2-planes in $\mathbb{R}^4$, there is a set of low-degree 2-surfaces and a set of low-degree 3-surfaces such that each surface has a degree less than $D$, and the surfaces in these sets contain all but $L^{\frac{3}{4}+\frac{1}{2}\epsilon}S + LS^{\frac{1}{2}+\epsilon}$ of the incidences between lines and 2-planes. 

\begin{theorem}
\label{improved}
Let there be a set of $L$ lines and a set of $S$ 2-planes in $\mathbb{R}^4$ where $L^{1/2} \ll S \ll L$, and a constant $D>1$. Let there also be a set of irreducible algebraic 2-surfaces, $\mathcal{X}$, so that each 2-surface in $\mathcal{X}$ has a degree of at most $D$ and contains at least $L^{\frac{1}{2}+\epsilon}$ lines, and a set of irreducible algebraic 3-surfaces, $\mathcal{Y}$, so that each 3-surface in $\mathcal{Y}$ has a degree of at most $D$ and contains at least $S^{\frac{1}{2}+\epsilon}$ 2-planes. Then, for an $\epsilon>0$, the number of intersections formed between lines that are not contained in the 2-surfaces in $\mathcal{X}$ and 2-planes that are not contained in the 3-surfaces in $\mathcal{Y}$ is at most $L^{\frac{3}{4}+\frac{1}{2}\epsilon}S + LS^{\frac{1}{2}+\epsilon}$.
\end{theorem}

Theorem \ref{improved} implies Theorem \ref{main}. The 2-surfaces in $\mathcal{X}$ and the 3-surfaces in $\mathcal{Y}$ violate the non-degeneracy conditions in Theorem \ref{main}. If both $\mathcal{X}$ and $\mathcal{Y}$ are empty, then the configuration of lines and planes is truly 4-dimensional, as there are no 2-surfaces containing more than $L^{1/2+\epsilon}$ lines and no 3-surfaces containing more than $S^{1/2+\epsilon}$ 2-planes. 

\section{Proof}
\subsection{Polynomial Partitioning}
Let $\mathcal{L}$ denote the set of $L$ lines and $\mathcal{S}$ denote the set of $S$ 2-planes in $\mathbb{R}^4$. We also have that $D>1$ and $\epsilon>0$. We start our proof for Theorem \ref{improved} by using polynomial partitioning. When partitioning, we have two choices in terms of how to choose the degree, $D$, of the partitioning polynomial. In earlier works by Guth and Katz when they developed the technique, $D$ was chosen to equal a power of the number of varieties in space. However, choosing $D$ in this fashion leads to a large value of $D$. This makes the zero-set of the partitioning polynomial very complicated, making controlling the intersections that lie on the polynomial surface difficult. However, this issue can be avoided by choosing $D$ to be equal to a constant, an approach developed by Solymosi and Tao in \cite{solymositao} when they generalized the Szemer\'edi-Trotter Theorem to higher dimensions. This approach was adopted by Guth and other authors in their recent works, and we will use it as well. Thus, by choosing $D$ to be a constant, we have that $D << L$ and $D << S$. We then partition the space by applying Guth's partitioning theorem for varieties twice: first with a polynomial $P_1$ of degree $D_1$ to divide lines evenly into cells, and then with a polynomial $P_2$ of degree $D_2$ to divide 2-planes evenly. We then multiply the two polynomials together to obtain a new polynomial $P=P_1P_2$ with degree $D=D_1D_2$ that partitions both the lines and the 2-planes evenly into cells. $P$ divides $\mathbb{R}^4$ into connected components of $\mathbb{R}^4 \setminus Z(P)$---where $Z(P)$ denotes the zero-set of $P$---resulting in roughly $D^4$ open cells. The exponent $4$ in $D^4$ corresponds to the dimension of space $\mathbb{R}^4$.

In order to use polynomial partitioning and induction, we need to estimate how many lines and 2-planes will be in each cell $O_i$. By the principal of polynomial partitioning technique, $P$ divides lines and 2-planes evenly among open cells. According to B\'ezout's Theorem, the number of intersections between a line and the polynomial surface with degree $D$ is bounded by $D$. Thus, one single line can enter at most $D+1$ cells, and the number of lines of $\mathcal{L}$ that intersect each open cell is bounded by $\frac{L(D+1)}{D^4}$, or approximately $\frac{L}{D^3}$. Similarly, to estimate the number of open cells a single 2-plane can enter, consider the intersection of a 2-plane with degree one and the partitioning polynomial surface of degree $D$: their intersection is a curve with degree $D$ (from B\'ezout's Theorem), and has at most $D^2$ self-intersections (from B\'ezout's Theorem), hence dividing the 2-plane into at most $D^2$ regions. Thus, each 2-plane can enter at most $D^2$ cells, and the number of 2-planes of $\mathcal{S}$ that intersect each open cell is bounded by $\frac{SD^2}{D^4}$, or $\frac{S}{D^2}$.

Now that the lines and 2-planes are divided into $D^4$ open cells approximately evenly, we count the intersections in open cells and the intersections on the zero-set surface of the partitioning polynomial separately, and then add them together.

\subsection{Intersections in Open Cells}
In this section we bound the intersections in the $D^4$ open cells. Let $\mathcal{X}$ be a set of 2-surfaces, each containing at least $L^{1/2+\epsilon}$ lines, and let $\mathcal{Y}$ be a set of 3-surfaces, each containing at least $S^{1/2+\epsilon}$ 2-planes. In order to use induction, we let $\mathcal{X}_i$ to be the set of 2-surfaces in each open cell $O_i$ where each 2-surface contains at least $L_i^{1/2+\epsilon} = (\frac{L}{D^3})^{1/2+\epsilon}$ lines, and let $\mathcal{Y}_i$ to be the set of 3-surfaces in each open cell $O_i$, where each 3-surface in $\mathcal{Y}_i$ contains at least $S_i^{1/2+\epsilon} = (\frac{S}{D^2})^{1/2+\epsilon}$ 2-planes.

After applying induction to each cell, the number of intersections in each cell $O_i$ formed between the lines that are not contained in the 2-surfaces in $\mathcal{X}_i$ and the 2-planes that are not contained in the 3-surfaces in $\mathcal{Y}_i$ is bounded by $L_i^{\frac{3}{4}+\frac{1}{2}\epsilon}S_i + L_iS_i^{\frac{1}{2}+\epsilon}$;
substituting $\frac{L}{D^3}$ for $L_i$ and $\frac{S}{D^2}$ for $S_i$ gives $\frac{L}{D^3}^{\frac{3}{4}+\frac{1}{2}\epsilon}\frac{S}{D^2} + \frac{L}{D^3}\frac{S}{D^2}^{\frac{1}{2}+\epsilon}$.

Next, we sum the intersections in all cells by multiplying the bound for an individual cell, $\frac{L}{D^3}^{\frac{3}{4}+\frac{1}{2}\epsilon}\frac{S}{D^2} + \frac{L}{D^3}\frac{S}{D^2}^{\frac{1}{2}+\epsilon}$, by the number of open cells, $D^4$, which gives $\frac{L}{D^3}^{\frac{3}{4}+\frac{1}{2}\epsilon}SD^2+LD\frac{S}{D^2}^{\frac{1}{2}+\epsilon}$; simplifying, the result is $D^{-1/4-3/2\epsilon}(L^{\frac{3}{4}+\frac{1}{2}\epsilon}S) + D^{-2\epsilon}(LS^{\frac{1}{2}+\epsilon})$. 

As $D$ is a large constant and $\epsilon$ is positive, the exponents of the $D$ terms are negative, and the above bound is better than our desired bound of $L^{\frac{3}{4}+\frac{1}{2}\epsilon}S + LS^{\frac{1}{2}+\epsilon}$.

\subsection{Pruning $\mathcal{X}$ and $\mathcal{Y}$} 
In the previous step, after summing up $\mathcal{X}_i$ and $\mathcal{Y}_i$ for all cells, we are left with the sets $\mathcal{X}$ and $\mathcal{Y}$ containing 2-surfaces and 3-surfaces from the entire space. The following conditions were used to define the 2-surfaces and 3-surfaces of $\mathcal{X}_i$ and $\mathcal{Y}_i$: each 2-surface in $\mathcal{X}_i$ contains at least ${L_i}^{1/2+\epsilon}=(\frac{L}{D^3})^{1/2+\epsilon}$ lines, and each 3-surface in $\mathcal{Y}_i$ contains at least ${S_i}^{1/2+\epsilon}=(\frac{S}{D^2})^{1/2+\epsilon}$ 2-planes. These conditions, however, are much smaller than the conditions given in Theorem \ref{improved} for $\mathcal{X}$ and $\mathcal{Y}$, which are: 
each 2-surface in $\mathcal{X}$ contains at least $L^{1/2+\epsilon}$ lines, and each 3-surface in $\mathcal{Y}$ contains at least $S^{1/2+\epsilon}$ 2-planes. That is, $(\frac{L}{D^3})^{1/2+\epsilon}$ is much smaller than $L^{1/2+\epsilon}$ and $(\frac{S}{D^2})^{1/2+\epsilon}$ is much smaller than $S^{1/2+\epsilon}$. We have included 2-surfaces and 3-surfaces in the construction of $\mathcal{X}_i$ and $\mathcal{Y}_i$ for each cell that do not satisfy the conditions in Theorem \ref{improved}. More specifically, after summing $\mathcal{X}_i$ for each cell to obtain $\mathcal{X}$, there are still some 2-surfaces $X_i$ that contain between $\frac{L}{D^3}^{1/2+\epsilon}$ and $L^{1/2+\epsilon}$ lines. These intersections contributed by these 2-surfaces have not been considered in the last section yet, but should be since they contain less than $L^{1/2+\epsilon}$ lines. Similarly, in $\mathcal{Y}$ there are still some 3-surfaces $Y_i$ that contain between $\frac{S}{D^2}^{1/2+\epsilon}$ and $S^{1/2+\epsilon}$ 2-planes which contain intersections that should be counted. In the following sections, we bound the intersections contained in these 2-surfaces and 3-surfaces and add their contributions to the total bound.

\subsection{Intersections with 2-Surfaces}
First, we determine the total number of 2-surfaces, denoted by $G_2$, that contain between $(\frac{L}{D^3})^{\frac{1}{2}+\epsilon}$ and $L^{\frac{1}{2}+\epsilon}$ lines. According to B\'ezout's Theorem, two 2-surfaces of degree $D$ can share up to $D^2$ lines; therefore, we cannot use $L/[(\frac{L}{D^3})^{\frac{1}{2}+\epsilon}]$ to bound $G_2$---it must be larger. Guth proved a lemma for this same question \cite{guthdistinctdistances}.

\begin{lemma}
\label{twolemma}
Suppose $\mathcal{L}$ is a set of lines in $\mathbb{R}^3$, and $\mathcal{X}$ is a set of irreducible algebraic surfaces of degree at most $D$, and suppose that each surface $X \in \mathcal{X}$ contains at least $A$ lines of $\mathcal{L}$. 
If $A>2D|\mathcal{L}|^{1/2}$, then $|\mathcal{X}| \leq \frac{2|\mathcal{L}|}{A}$
\end{lemma}

Using the same proof given by Guth in his paper, we can prove that this lemma also applies to $\mathbb{R}^4$. We can then apply this lemma to our 2-surface situation. As each 2-surface contains at least $\frac{L}{D^3}^{\frac{1}{2}+\epsilon}$ lines, we let $A=\frac{L}{D^3}^{\frac{1}{2}+\epsilon}$. We choose an $\epsilon$ so that $L^\epsilon > 2D^{\frac{5}{2}+3\epsilon}$. $D$ is a constant, $D \ll L$, and $\epsilon$ is positive, so this condition is met. Therefore $A>2DL^{\frac{1}{2}}$, and the lemma gives $G_2 \leq \frac{2L}{A} < \frac{2L}{(L/D^3)^{1/2+\epsilon}}$; which simplifies to $G_2 \leq \frac{2L}{A} < 2D^{\frac{3}{2}+3\epsilon}L^{\frac{1}{2}-\epsilon}$.

To bound the incidences contained in these 2-surfaces, it is necessary to consider all possible configurations: the intersection between a 2-surface and a 2-plane from $\mathcal{S}$ may be a set of points, a curve, or a plane.

\begin{case2surface} For 2-planes who intersect 2-surfaces at a finite number of points, the incidence between 2-planes and lines at the 2-surfaces is bounded by $2D^{\frac{5}{2}+3\epsilon}(L^{\frac{1}{2}-\epsilon}S)$.
\end{case2surface}

\begin{proof}
From B\'ezout's theorem, one 2-plane can intersect one 2-surface of degree $D$ at at most at $D$ points. There are at most $S$ 2-planes, and at most $2D^{\frac{3}{2}+3\epsilon}L^{\frac{1}{2}-\epsilon}$ 2-surfaces, so the bound is $D*S*(2D^{3/2+3\epsilon}L^{1/2-\epsilon})=2D^{\frac{5}{2}+3\epsilon}(L^{\frac{1}{2}-\epsilon}S)$.
\end{proof}

\begin{case2surface}
For 2-planes whose intersections with 2-surfaces form a curve, the incidence bound is $(C_3D^{\frac{3}{2}+3\epsilon})(LS^\frac{1}{2})$, where $C_3$ is a constant depending on $\epsilon$.
\end{case2surface}

\begin{proof}
The proof of this bound is somewhat involved. We need to figure out how many curves will be formed by intersections between $S$ 2-planes and $G_2$ 2-surfaces. In \cite{guthzahl}, Guth and Zahl proved the following theorem for a similar situation for 2-surfaces in $\mathbb{R}^4$. 

\begin{theorem}
For each $D \geq 1$, there are constants $C_1$, $C_2$ so that the following holds. Let $\mathcal{N}$ be a collection of $n$ irreducible 2-surfaces in $\mathbb{R}^4$, each of degree at most $D$, and let $A \geq C_1n^\frac{1}{2}$. Suppose that for each $N \in \mathcal{N}$, there are at least $C_2A$ distinct irreducible curves $\gamma \subset N$ that are incident to at least one other surface from $\mathcal{N}$. Then there is an irreducible three dimensional variety $M$ of degree at most $100D^2$ that contains $\geq A$ 2-surfaces from $\mathcal{N}$.
\end{theorem}
Essentially, this theorem says that if every 2-surface has at least $C_2A$ curves which are incident to at least one other 2-surface from $\mathcal{N}$, then there exists a 3-surface $M$ that contains at least $A(\geq C_1n^\frac{1}{2})$ 2-surfaces from $\mathcal{N}$. As there are $n$ 2-surfaces, the total number of these "2-rich curves", curves formed by the intersections of at least two 2-surfaces, will be $\frac{C_2An}{2} \geq \frac{C_2(C_1n^\frac{1}{2})n}{2}=\frac{C_1C_2}{2}n^\frac{3}{2}$, where the divisor 2 is used to account for the fact that two 2-surfaces count the same 2-rich curve twice. As the contrapositive of this statement is also true, we can state it as follows.
\begin{lemma}
\label{contra}
In above theorem, if no 3-hypersurface of degree at most $100D^2$ contains more than $C_1n^\frac{1}{2}$ 2-surfaces from $\mathcal{N}$ in $\mathbb{R}^4$, then the total number of 2-rich curves formed by the intersections between 2-surfaces from $\mathcal{N}$ has an upper bound of $\frac{C_1C_2}{2}n^\frac{3}{2}$. 
\end{lemma} 
We cannot use this lemma directly, because we are bounding the incidences between two separate sets of 2-dimensional objects (2-surfaces and 2-planes), not between objects from the same set of 2-surfaces. To count the intersections in our case, we need to develop a bipartite version between 2 sets of objects: $S$ 2-planes and $G_2$ 2-surfaces.

We combine those $G_2$ 2-surfaces and $S$ 2-planes to form a new set of $n$ 2-dimensional objects, $n=S+G_2$. We choose constant $C_1$ and $\epsilon$ so that $ S^{2\epsilon} < C_1^2 $, which gives $(S^{2\epsilon})S < C_1^2 S < C_1^2 S + C_1^2 G_2=C_1^2(S + G_2)=C_1^2n$. Taking the square root, we get $S^{1/2+\epsilon} < C_1n^{1/2}$. With this inequality, our non-degeneracy condition "no 3-surface contains more than $S^{\frac{1}{2}+\epsilon}$ 2-planes" will become "no 3-surface contains more than $C_1n^{1/2}$ 2-planes". The degree of 3-surfaces in our case is at most $D$, less than the requirement of at most $100D^2$ degree in Lemma \ref{contra}. Therefore both hypotheses of Lemma \ref{contra} are satisfied, and we can apply the lemma to our combined set of $n=S+G_2$ objects, giving a bound $I(n)$ for $n$: $I(n)=\frac{C_1C_2}{2}n^\frac{3}{2}=\frac{C_1C_2}{2}(S+G_2)^\frac{3}{2}$. To expand the exponent term, we can use the binomial theorem: 
$$(x+y)^r=\sum_{i=0}^{\infty}\frac{r(r-1)(r-2)....(r-i+1)}{i!}x^{r-i}y^i$$
$$I(n)=\frac{C_1C_2}{2}(S+G_2)^{\frac{3}{2}}=\frac{C_1C_2}{2}S^\frac{3}{2}(1+\frac{G_2}{S})^\frac{3}{2}=\frac{C_1C_2}{2}S^{\frac{3}{2}}[1+\frac{3}{2}(\frac{G_2}{S})+\frac{3}{8}(\frac{G_2}{S})^2-\frac{3}{48}(\frac{G_2}{S})^3+....]$$
Since $G_2 \leq 2D^{\frac{3}{2}+3\epsilon}L^{\frac{1}{2}-\epsilon}$, and we have the condition $S \gg L^\frac{{1}}{2}$, so we have $S \gg G_2$ and $\frac{G_2}{S} \ll 1$. Thus the expression is dominated by the first two terms, and we can drop the rest of the terms on the right-hand side in our estimate, resulting in the following bound:
$$I(n) \eqsim \frac{C_1C_2}{2}S^{\frac{3}{2}}[1+\frac{3}{2}(\frac{G_2}{S}) ]= \frac{C_1C_2}{2}S^{\frac{3}{2}}+\frac{3C_1C_2}{4}S^\frac{1}{2}G_2=I_{plane} + I_{surface}+I_{plane-surface}$$

The above bound includes 3 groups of intersections: between two 2-planes from the set of $S$ 2-planes, $I_{plane}$; between two 2-surfaces from the set of $G_2$ 2-surfaces, $I_{surface}$; and between one 2-plane from $S$ 2-planes and one 2-surface from $G_2$ 2-surfaces, $I_{plane-surface}$. We only need the last term. $I_{plane}$ can be obtained by applying the lemma directly to the set of $S$ 2-planes alone, and its bound is simply $I_{plane}=\frac{C_1C_2}{2}S^{\frac{3}{2}}$, which equals the first term in above equation. Subtracting this term, and moving $I_{surface}$ to the other side of equation, we get:
$$ I_{plane-surface}  = \frac{3C_1C_2}{4}S^\frac{1}{2}G_2 - I_{surface} < \frac{3C_1C_2}{4}S^\frac{1}{2}G_2 = (\frac{3C_1C_2}{4}S^\frac{1}{2})(2D^{\frac{3}{2}+3\epsilon}L^{\frac{1}{2}-\epsilon})$$
Due to the $I_{surface}$ term in the inequality, this bound is not sharp and we can safely use this bound. Next, we multiply the bound by the maximum number of lines each 2-surface can contain, which is at most $L^{1/2+\epsilon}$. This gives us the bound: 
$$(\frac{3C_1C_2}{4}S^\frac{1}{2})(2D^{\frac{3}{2}+3\epsilon}L^{\frac{1}{2}-\epsilon})(L^{\frac{1}{2}+\epsilon})=(\frac{3C_1C_2}{2}D^{\frac{3}{2}+3\epsilon})(LS^\frac{1}{2})=(C_3D^{\frac{3}{2}+3\epsilon})(LS^\frac{1}{2})$$
where we define $C_3=3C_1C_2/2$.
\end{proof}

\begin{case2surface}
For 2-planes whose intersection with 2-surfaces is a plane, the number of intersections between 2-planes and lines is bounded by $(2D^{\frac{5}{2}+3\epsilon})L$.
\end{case2surface}

\begin{proof}
In this case, the 2-planes are completely contained within a 2-surface. By B\'ezout's Theorem, each 2-surface of degree $D$ can contain at most $D$ 2-planes. Each of the $2D^{\frac{3}{2}+3\epsilon}L^{\frac{1}{2}-\epsilon}$ 2-surfaces can contain at most $L^{1/2+\epsilon}$ lines, so the number of incidences is $D(2D^{\frac{3}{2}+3\epsilon}L^{\frac{1}{2}-\epsilon})(L^{\frac{1}{2}+\epsilon})$, which simplifies to $2D^{\frac{5}{2}+3\epsilon}L$. 
\end{proof}

\subsection{Intersections with 3-Surfaces}
As with 2-surfaces, we need to bound the maximum number of 3-surfaces, denoted by $G_3$, that contain between $\frac{S}{D^2}^{\frac{1}{2}+\epsilon}$ and $S^{\frac{1}{2}+\epsilon}$ 2-planes. In this case, we have following lemma for 3-surfaces, which is similar to Lemma \ref{twolemma} from \cite{guthdistinctdistances} for 2-surfaces.

\begin{lemma}
\label{threelemma}
Given a set of $S$ 2-planes in $\mathbb{R}^4$, and a set of $G_3$ irreducible algebraic 3-surfaces of degree at most $D$, each containing at least $A$ 2-planes: if $A>2DS^{1/2}$, then $G_3 \leq \frac{2S}{A}$.
\end{lemma}

Lemma \ref{threelemma} is almost identical to Lemma \ref{twolemma}, with $L$ lines being replaced by $S$ 2-planes and $G_2$ 2-surfaces being replaced by $G_3$ 3-surfaces. The maximum number of 2-planes that two 3-surfaces may share is also $D^2$, the same as the maximum number of lines that may be shared by two 2-surfaces. Thus, the proof is the same as the proof for Lemma \ref{twolemma}.

Lemma \ref{threelemma} gives us a bound for the number of 3-suraces. As each 3-surface contains at least $\frac{S}{D^2}^{\frac{1}{2}+\epsilon}$ 2-planes, we let $A=\frac{S}{D^2}^{\frac{1}{2}+\epsilon}$. An $\epsilon$ can be chosen so that $S^\epsilon > 2D^{2+2\epsilon}$. $D$ is a constant, $D \ll S$, and $\epsilon$ is positive, and therefore $A>2DS^{1/2}$. Thus, the lemma gives an upper bound for $G_3$ of $G_3 \leq \frac{2S}{A}$; substituting $A$ with $\frac{S}{D^2}^{\frac{1}{2}+\epsilon}$ and simplifying, the result is 
$G_3 \leq 2D^{1+2\epsilon}S^{\frac{1}{2}-\epsilon}$.

To bound the incidences contained in these 3-surfaces, we need to consider two situations.

\begin{case3surface}
For lines that are not contained in 3-surfaces, the number of intersections between those lines and the 2-planes contained in 3-surfaces is $(2D^{2+2\epsilon})(LS^{1/2-\epsilon})$.
\end{case3surface}

\begin{proof}
By B\'ezout's theorem, each line may intersect one 3-surface at at most $D$ points. As there are at most $L$ lines, and at most $2D^{1+2\epsilon}S^{\frac{1}{2}-\epsilon}$ 3-surfaces, the maximum number of intersections is $DL(2D^{1+2\epsilon}S^{1/2-\epsilon})$, which simplifies to $(2D^{2+2\epsilon})(LS^{1/2-\epsilon})$.
\end{proof}

\begin{case3surface}
For lines that are contained in 3-surfaces, the intersections between those lines and 2-planes contained in 3-surfaces is $(2D^{1+2\epsilon})(L^{\frac{3}{4}+\frac{1}{2}\epsilon}S) + (LS^{\frac{1}{2}+\epsilon})$. 
\end{case3surface}

\begin{proof}
To obtain this bound, we construct a bipartite graph consisting of two sets of vertices: one set of vertices $\mathcal{U}$ representing the $L$ lines, and the other set of vertices $\mathcal{V}$ representing the $G_3$ 3-surfaces. Each edge in this graph corresponds to a line, 3-surface pair such that the line is contained in the 3-surface. 

Their intersection is a 2-surface. If this 2-surface contains more than $L^{1/2+\epsilon}$ lines, then we will not count its intersections as we only count intersections in 2-surfaces containing no more than $L^{1/2+\epsilon}$ lines. This can be modeled by a complete bipartite subgraph of the form $K_{L^{1/2+\epsilon},2}$, where $L^{1/2+\epsilon}$ represents $L^{1/2+\epsilon}$ lines and $2$ represents two 3-surfaces. Because this bipartite graph is complete, and each of the $L^{1/2+\epsilon}$ vertices are connected to each of the $2$ vertices, this represents two 3-surfaces each containing the same $L^{1/2+\epsilon}$ lines. Therefore, the bipartite graph cannot contain any complete subgraphs of the form $K_{L^{1/2+\epsilon},2}$. 

The Kov\'ari-S\'os-Tur\'an Theorem \cite{kst} bounds the maximum number of edges of a bipartite graph, with two sets of vertices $\mathcal{U}$ and $\mathcal{V}$ for which $|\mathcal{U}|=m$ and $|\mathcal{V}|=n$, that does not contain any complete subgraph of the form $K_{s,t}$. The theorem is this: $z(m, n; s,t) \leq (s-1)^\frac{1}{t}(n-t+1)m^{1-\frac{1}{t}}+(t-1)m$. Substituting $m=L$, $n=G_3$, $s=(L^{1/2+\epsilon}+1)$, and $t=2$, the inequality becomes $z(L, G_3; (L^{\frac{1}{2} + \epsilon}+1), 2) \leq [(L^{\frac{1}{2}+\epsilon}+1)-1]^\frac{1}{2}(G_3-2+1)L^{1-\frac{1}{2}} + (2-1)L$. Simplifying, $z(L, G_3; (L^{\frac{1}{2} + \epsilon}+1), 2) \leq L^{\frac{1}{4} + \frac{1}{2}\epsilon}G_3L^\frac{1}{2}+L=L^{\frac{3}{4}+\frac{1}{2}\epsilon}G_3+L$. We multiply this bound by the maximum number of 2-planes, $S^{1/2+\epsilon}$, which each 3-surface may contain. This gives the bound $(L^{\frac{3}{4}+\frac{1}{2}\epsilon}G_3+L)(S^{\frac{1}{2}+\epsilon})\leq [L^{\frac{3}{4}+\frac{1}{2}\epsilon}(2D^{1+2\epsilon}S^{\frac{1}{2}-\epsilon})+L](S^{\frac{1}{2}+\epsilon}) = (2D^{1+2\epsilon})(L^{\frac{3}{4}+\frac{1}{2}\epsilon}S)+LS^{\frac{1}{2}+\epsilon}$.
\end{proof}

\subsection{Contributions from 2-Surfaces and 3-Surfaces}

Summing up contributions from all 2-surface cases and all 3-surface cases, we have:
$$I(L,S)_{2-surfaces}+I(L,S)_{3-surfaces}=2D^{\frac{5}{2}+3\epsilon}(L^{\frac{1}{2}-\epsilon}S)+(C_3D^{\frac{3}{2}+3\epsilon})(LS^\frac{1}{2})$$
$$+(2D^{\frac{5}{2}+3\epsilon})L+(2D^{2+2\epsilon})(LS^{1/2-\epsilon})+(2D^{1+2\epsilon})(L^{3/4+1/2\epsilon}S) + (LS^{1/2+\epsilon})$$
$$=O(L^{3/4+1/2\epsilon}S + LS^{1/2+\epsilon})$$
which is the bound we started our induction with.

\subsection{Intersections on the Zero-Set of the Partitioning Polynomial}
Lastly, we bound the intersections at the zero-set of the partitioning polynomial $Z(P)$. In some cases there will be lines that intersect, but are not contained in the zero-set; in other cases the lines will be contained in the zero-set, and then it depends on how the 2-planes intersect the zero-set to bound the line-plane incidences.

\begin{casezero}
For lines that are not contained in the zero-set, each line may intersect the zero-set of degree $D$ at most $D$ times by B\'ezout's Theorem, giving at most $D$ intersections. Thus, the bound is $DL$.
\end{casezero}

\begin{casezero} 
For lines that are completely contained within the zero-set and 2-planes whose intersection with the zero-set is a set of points, each 2-plane can be tangent to the zero-set at at most at $D$ points by B\'ezout's Theorem. Thus, the bound is $DS$.
\end{casezero}

\begin{casezero} 
For lines that are completely contained within the zero-set and 2-planes whose intersection with the zero-set form a curve (1-dimensional case), the number of intersections between these lines and 2-planes is at most $(\frac{3}{8}+\frac{\epsilon}{4})C_4(L^{\frac{1}{2}+\epsilon}S)$, where $C_4$ is a constant depending on $\epsilon$
\end{casezero}

\begin{proof}
As there are $S$ 2-planes, the number of curves, $V$, is at most $S$, and $ V \leq S$. The degree of these curves is at most $D$, same as the degree of zero-set, as per B\'ezout's Theorem. We need to count the intersections between $L$ lines and these $V$ curves. In \cite{guthdistinctdistances}, Guth proved a similar bound.

\begin{theorem}
For any $\epsilon>0$, there are $D(\epsilon), C_4(\epsilon)$ so that the following holds. If $\mathcal{N}$ is a set of $n$ lines in $\mathbb{R}^3$, and there are less than $n^{\frac{1}{2}+\epsilon}$ lines in any algebraic surface of degree at most D, and if $2 \leq r \leq 2n^\frac{1}{2}$, the number of r-rich points is bounded by $C_4n^{\frac{3}{2}+\epsilon}r^{-2}$. (Note: when $r=2$, the bound becomes: the number of 2-rich points is bounded by $\frac{C_4}{4}n^{\frac{3}{2}+\epsilon}$).
\end{theorem}

Here, "r-rich point" denotes an intersection point where r lines intersect. The above bound is for $\mathbb{R}^3$. In our case, both 1-dimensional objects (lines and curves) are contained within the zero-set $Z(P)$, which by itself is a 3-dimensional sub-space; so our case is similar to the above theorem. This theorem has also been generalized to curves in \cite{guthzahl} and \cite{guthbook}. Nevertheless, we cannot apply the theorem to our case directly, and need to develop a bipartite version, because we want to count the intersections between two separate sets of objects (between $L$ lines and $V$ curves), not between two objects from the same set of objects.

Again, we use the same process as we did in Case 2 of Intersections with 2-Surfaces by combining $L$ lines and $V$ curves to form a new set of $n=L+V$ objects. Our non-degeneracy condition says that there are less than $L^{1/2+\epsilon}$ lines in any 2-surface, so it is also true that there are less than $n^{1/2+\epsilon}$ lines in any 2-surface as $L<n$; thus, the conditions of the theorem is satisfied. We apply the theorem and get an upper bound of 2-rich points for $n=L+V$; $I=\frac{C_4}{4}n^{\frac{3}{2}+\epsilon}=\frac{C_4}{4}(L+V)^{\frac{3}{2}+\epsilon}$. As we have the condition $L \gg S$ in Theorem \ref{main}, and $S \geq V$, we get $L \gg V$ and $L/V \ll 1$. Then we can expand this using the binomial theorem to get
$$I=\frac{C_4}{4}(L+V)^{\frac{3}{2}+\epsilon}\simeq \frac{C_4}{4}L^{\frac{3}{2}+\epsilon}+\frac{C_4}{4}(\frac{3}{2}+\epsilon)L^{\frac{1}{2}+\epsilon}V = I_{line}+I_{curve}+I_{line-curve}$$
Separately, if we apply the theorem to the set of $L$ lines alone, we get the 2-rich point bound for the set of $L$ lines alone as $I_{line}=\frac{C_4}{4}L^{\frac{3}{2}+\epsilon}$. Subtracting $I_{line}$ from this equation, we obtain our desired bound for 2-rich points between line and curve, $I_{line-curve}$, as
$$I_{line-curve}= \frac{C_4}{4}(\frac{3}{2}+\epsilon)L^{\frac{1}{2}+\epsilon}V - I_{curve} < (\frac{3}{8}+\frac{\epsilon}{4})C_4L^{\frac{1}{2}+\epsilon}V \leq (\frac{3}{8}+\frac{\epsilon}{4})C_4L^{\frac{1}{2}+\epsilon}S$$ 
The term $I_{curve}$, accounting for the bound for the set of $V$ alone, is ignored in above inequality. Therefore the above bound of $I_{line-curve}$ is not sharp, and we can safely use this bound in our proof. As a curve is a subset of a 2-plane, it can intersect with a line at most once. Thus, the above inequality is our final bound for this case.
\end{proof}
\begin{casezero}
For lines that are completely contained within the zero-set and 2-planes whose intersection with the zero-set form a plane: in this case the 2-planes must be completely contained within the zero-set, and the number of intersections between these lines and 2-planes is at most $LS^{1/2+\epsilon}$.
\end{casezero}

\begin{proof}
The zero-set itself is a 3-surface. The non-degeneracy condition says that any 3-surface can contain at most $S^{1/2+\epsilon}$ 2-planes. Therefore, there are at most $S^{\frac{1}{2}+\epsilon}$ 2-planes in this case. As there are at most $L$ lines, the number of intersections between lines and 2-planes that are completely contained within the zero-set is bounded by $LS^{1/2+\epsilon}$.
\end{proof}

\subsubsection{Summing Cases of the Zero-Set}
Combining intersection contributions from the above four cases, the total number of intersections on the zero-set of the partitioning polynomial is at most
$$I(L,S)_{zero-set}= DL + DS + [(\frac{3}{8}+\frac{\epsilon}{4})C_4L^{\frac{1}{2}+\epsilon}S] + LS^{\frac{1}{2} + \epsilon}  \simeq  (\frac{3}{8}+\frac{\epsilon}{4})C_4L^{\frac{1}{2}+\epsilon}S + LS^{\frac{1}{2} + \epsilon}$$ 
which is considerably smaller than our bound in open cells and in 2-surfaces and 3-surfaces. This means the contribution from zero-set is relatively little compared to other sources.

\subsection{Total Incidences}
Adding up all the contributions from the open cells, the pruning of $\textit{X}$ and $\textit{Y}$, and the zero-set, we have the final bound
$$I(L,S)_{cells}+I(L,S)_{surfaces}+I(L,S)_{zero-set}$$
$$=D^{-1/4-3/2\epsilon}(L^{3/4+1/2\epsilon}S) + D^{-2\epsilon}(LS^{1/2+\epsilon})+(L^{3/4+1/2\epsilon}S + LS^{1/2+\epsilon})+(\frac{3}{8}+\frac{\epsilon}{4})C_4L^{\frac{1}{2}+\epsilon}S + LS^{\frac{1}{2} + \epsilon}$$
$$=O{(L^{\frac{3}{4}+\frac{1}{2}\epsilon}S + LS^{\frac{1}{2}+\epsilon})}$$ 

which is the bound we started our induction with. This closes the induction and proves Theorem \ref{improved}. When $\mathcal{X}$ and $\mathcal{Y}$ are empty, which signifies that both non-degeneracy conditions are met, Theorem \ref{improved} is reduced to the main Theorem \ref{main}.

\section{Conclusion}
Having considered incidences within cells after polynomial partitioning, incidences in 2-surfaces and 3-surfaces containing large numbers of lines and planes, and incidences contained in the zero-set of the partitioning polynomial, the total contribution from all these cases is less than or equal to our desired bound of $O(L^{\frac{3}{4}+\frac{1}{2}\epsilon}S + LS^{\frac{1}{2}+\epsilon})$, and thus we have proven that for an $\epsilon>0$, $L^{\frac{3}{4}+\frac{1}{2}\epsilon}S + LS^{\frac{1}{2}+\epsilon}$ is the upper bound for the number of intersections between $L$ lines and $S$ planes in $\mathbb{R}^4$.

\section{Future Directions}
Although this paper finds a non-trivial bound for the intersections between lines and 2-planes in $\mathbb{R}^4$, this is merely the first step in studying line-plane incidences in $\mathbb{R}^4$. The most obvious future direction is to find a sharper bound; the $\epsilon$ that produces the sharpest bound has not yet been determined. This may be accomplished by exploring other methods of categorically counting incidences. Other explorations present themselves by changing the objects in the question. A bound for the intersections between lines and 3-planes in $\mathbb{R}^5$, for instance, is a novel idea. Another similar problem, by expanding Guth's theorem on 2-rich points on curves in $\mathbb{R}^4$, is bounding the intersections between 2-planes in $\mathbb{R}^5$.  Finally, in this paper we only consider shapes with a degree of one, such as lines or 2-planes. It is worthwhile to consider how higher degree shapes such as curves or 2-surfaces might behave in similar situations.

\thispagestyle{empty}
\end{document}